\documentclass{article}
\usepackage{caption,subcaption}
\pdfoutput=1 
\usepackage{amsmath, amsthm, amssymb, amstext, amsfonts}
\usepackage{enumerate}
\usepackage[applemac]{inputenc}
\usepackage{graphicx}

\theoremstyle{plain}
\newtheorem{thm}{Theorem}[section]
\newtheorem{lem}[thm]{Lemma}

\newtheorem{cor}[thm]{Corollary}

\theoremstyle{definition}
\newtheorem{defi}[thm]{Definition}
\newtheorem{ex}{Example}

\newcount\commentno
\def\COMMENT#1{$^{<\the\commentno>}$%
     \vadjust{\vbox to 0pt{\vss\vskip-8pt\rightline{%
     \rlap{\hbox{\hskip7mm \vbox{\pretolerance=-1
     \doublehyphendemerits=0 \finalhyphendemerits=0
     \hsize40mm\tolerance=10000\eightpoint
     \lineskip=0pt\lineskiplimit=0pt
     \rightskip=0pt plus16mm\baselineskip8pt\noindent
     \hskip0pt       
     {$\langle$\the\commentno. #1$\rangle$}\endgraf}}}}\vss}}%
     \global\advance\commentno by1}%
\def\writecommentsasfootnotes{%
 \def\COMMENT{\global\advance\commentno by1\footnote{$^{<\the\commentno>}$}}%
 }
\def\nocomments{\def\COMMENT##1{}}
\def\?#1{\vadjust{\vbox to 0pt{\vss\vskip-8pt\leftline{%
     \llap{\hbox{\vbox{\pretolerance=-1
     \doublehyphendemerits=0\finalhyphendemerits=0
     \hsize16truemm\tolerance=10000\small
     \lineskip=0pt\lineskiplimit=0pt
     \rightskip=0pt plus16truemm\baselineskip8pt\noindent
     \hskip0pt        
     #1\endgraf}\hskip7truemm}}}\vss}}}
\newenvironment{txteq}
  {
    \begin{equation}
    \begin{minipage}[c]{0.85\textwidth} 
    \em                                
  }
  {\end{minipage}\end{equation}\ignorespacesafterend}
\newenvironment{txteq*}
  {
    \begin{equation*}
    \begin{minipage}[c]{0.85\textwidth} 
    \em                                
  }
  {\end{minipage}\end{equation*}\ignorespacesafterend}
\def\specrel#1#2{\mathrel{\mathop{\kern0pt #1}\limits_{#2}}}
\def\Specrel#1#2{\mathrel{\mathop{\kern0pt #1}\limits^{#2}}}

\newcommand{\lk}{$({<k)}$}

\newcommand{\N}{\ensuremath{\mathbb{N}}}

\newcommand{\Lra}{\ensuremath{\Leftrightarrow}}

\newcommand{\sm}{\ensuremath{\smallsetminus}}

\newcommand{\es}{\ensuremath{\emptyset}}

\newcommand{\sub}{\subseteq}

\newcommand{\cN}{{\ensuremath N}}
\newcommand{\cO}{\ensuremath{\mathcal O}}
\newcommand{\cP}{\ensuremath{\mathcal P}}
\newcommand{\cQ}{\ensuremath{\mathcal Q}}
\newcommand{\cR}{{\ensuremath R}}
\newcommand{\cS}{{\ensuremath S}}
\newcommand{\cT}{{\ensuremath T}}

\newcommand{\cV}{\ensuremath{\mathcal V}}

\newcommand{\ops}{\ensuremath{\sext, \sloc, \smax, \srext, \srloc, \srmax}}
\newcommand{\rops}{\ensuremath{\srext, \srloc, \srmax}}
\newcommand{\nops}{\ensuremath{\sext, \sloc, \smax}}

\newcommand{\sfont}[1]{\ensuremath{\mathsf{#1}}}

\newcommand{\srmax}{\sfont{all_r}}
\newcommand{\srloc}{\sfont{loc_r}}
\newcommand{\srext}{\sfont{ext_r}}

\newcommand{\smax}{\sfont{all}}
\newcommand{\sloc}{\sfont{loc}}
\newcommand{\sext}{\sfont{ext}}

\newcommand{\sMax}{\sfont{All}}
\newcommand{\sLoc}{\sfont{Loc}}
\newcommand{\sExt}{\sfont{Ext}}

\newcommand{\SP}{\ensuremath{(\cS,\cP)}}

\newcommand{\RP}{\ensuremath{(\cR,\cP)}}

\newcommand{\SPi}[1]{\ensuremath{(\cS_{#1},\cP_{#1})}}

\newcommand{\TV}{\ensuremath{(\cT,\cV)}}

\def\td{tree-decom\-po\-si\-tion}

\newcommand{\sys}{separation system}

\newcommand{\sepn}[2]{\ensuremath{{(#1,#2)}}}
\newcommand{\AB}{\sepn AB}
\newcommand{\BA}{\sepn BA}
\newcommand{\CD}{\sepn CD}
\newcommand{\DC}{\sepn DC}
\newcommand{\EF}{\sepn EF}
\newcommand{\FE}{\sepn FE}
\newcommand{\XY}{\sepn XY}
\newcommand{\YX}{\sepn YX}

\let\doublebar=\|

\newcommand{\mcm}[3]{\newcommand{#1}[#2]{{\ensuremath{#3}}}} 

\mcm{\tuple}{1}{\langle #1 \rangle}
\mcm{\name}{1}{\ulcorner #1 \urcorner}
\mcm{\Nbb}{0}{\mathbb{N}}
\mcm{\Zbb}{0}{\mathbb{Z}}
\mcm{\Rbb}{0}{\mathbb{R}}
\mcm{\Cbb}{0}{\mathbb{C}}
\mcm{\Fbb}{0}{\mathbb{F}}
\mcm{\Bcal}{0}{\cal B}
\mcm{\Ccal}{0}{\cal C}
\mcm{\Dcal}{0}{\cal D}
\mcm{\Ecal}{0}{\cal E}
\mcm{\Fcal}{0}{\cal F}
\mcm{\Gcal}{0}{\cal G}
\mcm{\Hcal}{0}{\cal H}
\mcm{\Ical}{0}{\cal I}
\mcm{\Lcal}{0}{\cal L}
\mcm{\Mcal}{0}{\cal M}
\mcm{\Ncal}{0}{\cal N\!}
\mcm{\Ocal}{0}{\cal O}
\mcm{\Pcal}{0}{{\cal P}}
\mcm{\Scal}{0}{{\cal S}}
\mcm{\Tcal}{0}{{\cal T}}
\mcm{\Ucal}{0}{{\cal U}}
\mcm{\Vcal}{0}{{\cal V}}
\mcm{\Wcal}{0}{{\cal W}}
\mcm{\Ycal}{0}{{\cal Y}}
\mcm{\Mfrak}{0}{\mathfrak M}

\usepackage{verbatim}

\nocomments

\title{Canonical \td s of finite graphs\\ I. Existence and algorithms}
\author{J.\ Carmesin \and R.\ Diestel \and M.\ Hamann \and F.\ Hundertmark}

\begin{document}

\maketitle

\begin{abstract}
We construct \td s of graphs that distinguish all their $k$-blocks and tangles of order~$k$, for any fixed integer~$k$. We describe a family of algorithms to construct such decompositions, seeking to maximize their diversity subject to the requirement that they commute with graph isomorphisms. In particular, all the decompositions constructed are invariant under the automorphisms of the graph.
\end{abstract}

\section{Introduction}
A~\emph{$k$-block} in a graph $G$, where $k$ is any positive integer, is a maximal set $X$ of at least~$k$ vertices such that no two vertices $x,x'\in X$ can be separated in~$G$ by fewer than~$k$ vertices other than $x$ and~$x'$. Thus, $k$-blocks can be thought of as highly connected pieces of a graph, but their connectivity is measured not in the subgraph they induce but in the ambient graph.

Extending results of Tutte~\cite{TutteGrTh} and of Dunwoody and Kr\"on~\cite{DunwoodyKroenArXiv}, three of us and Maya Stein showed that every finite graph $G$ admits, for every integer~$k$, a \td\ $\TV$ of adhesion~$<k$ that distinguishes all its $k$-blocks~\cite{confing}. These decompositions are {\em canonical\/} in that the map ${G\mapsto \TV}$ commutes with graph isomorphisms.%
   \COMMENT{}
   In particular, the decomposition $\TV$ constructed for~$G$ is invariant under the automorphisms of~$G$.

Our next aim, then, was to find out more about the \td s whose existence we had just proved. What can we say about their parts? Will every part contain a $k$-block? Will those that do consist of just their $k$-block, or might they also contain some `junk'? Such questions are not only natural; their answers will also have an impact on the extent to which our tree-decompositions can be used for an obvious potential application, to the graph isomorphism prob\-lem in complexity theory. See Grohe and Marx~\cite{GroheMarxTKn} for recent progress on this.

When we analysed our existence proof in view of these questions, we found that even within the strict limitations imposed by canonicity we can make choices that will have an impact on the answers. For example, we can obtain different decompositions (all canonical) if we seek to, alternatively,  minimize the number of inessential parts,  minimize the sizes of the parts,\vadjust{\penalty-200} or just of the essential parts, or achieve a reasonable balance between these properties. (A~part is called \emph{essential} if it contains a $k$-block, and \emph{inessential} otherwise.)

In this paper we describe a large family of algorithms%
   \footnote{We should point out that our reason for thinking in terms of algorithms is not, at this stage, one of complexity considerations: these are interesting, but they are not our focus here. Describing a decomposition in terms of the algorithm that produces it is simply the most intuitive way to ensure that it will be canonical: as long as the instructions of how to obtain the decomposition refer only to invariants of the graph (rather than, say, to a vertex enumeration that has to be chosen arbitrarily at some point), the decomposition that this algorithm produces will also be an invariant.}
   that each produce a canonical \td\ for given $G$ and~$k$. Their parameters can be tuned to optimize this \td\ in terms of criteria such as those above.%
   \COMMENT{}
   In~\cite{CDHH13CanonicalParts} we shall apply these results to specify algorithms from the family described here for which we can give sharp bounds on the number of inessential parts, or which under specified conditions ensure that some or all essential parts consist only of the corresponding $k$-block.%
   \COMMENT{}

The existence theorems which our algorithms imply will extend our results from~\cite{confing} in that the decompositions constructed will not only distinguish all the $k$-blocks of a graph, but also its tangles of order~$k$. (Tangles were introduced by Robertson and Seymour~\cite{GMX} and can also be thought of as indicating highly connected parts of a graph.) In order to treat blocks and tangles in a unified way, we work with a common generalization called `profiles'. These appear to be of interest in their own right, as a way of locating desirable local substructures in very general discrete structures. More about profiles, including further generalizations of our existence theorems to such general structures (including matroids), can be found in~\cite{profiles}. More on $k$-blocks, including different kinds of examples, their relationship to tangles, the algorithmic complexity of finding them, and a block-decomposition duality theorem, can be found in~\cite{ForcingBlocks}.

All graphs in this paper will be finite, undirected and simple. Any  graph-theoretic terms not defined here are explained in~\cite{DiestelBook10noEE}. Unless otherwise mentioned, $G$ will denote an arbitrary finite graph with vertex set~$V$.

\section{Separation systems}\label{sec_nsys}
A pair $(A,B)$ of subsets of $V$ such that $A\cup B = V$ is called a \emph{separation} of~$G$ if there is no edge $e = \{x,y\}$ in $G$ with $x \in A\sm B$ and $y\in B\sm A$. A separation $(A,B)$ is \emph{proper} if neither $A \sub B$ nor $B\sub A$; otherwise it is \emph{improper}. The \emph{order} of a separation \AB\ is the cardinality of its \emph{separator} $A\cap B$. A~separation of order $k$ is a \emph{$k$-separation}. By a simple calculation we obtain:

\begin{lem}\label{counting}
For any two separations $(A,B)$ and $(C,D)$, the orders of the separations $(A\cap C, B\cup D)$ and $(B\cap D, A\cup C)$ sum to $|A\cap B|+|C\cap D|$.\qed
\end{lem}

\goodbreak

We define a partial ordering on the set of separations of $G$ by
\begin{equation}\label{eq_order}
\AB\le\CD :\Lra A\sub C \wedge B \supseteq D.
\end{equation}
A separation \AB\ is \emph{nested} with a separation \CD, written as $\AB\|\CD$, if it is $\le$-comparable with either \CD\ or \DC. Since
\begin{equation}\label{eq_sym}
\AB\le\CD \Lra \DC\le\BA,
\end{equation}
the relation $\parallel$ is reflexive and symmetric.\footnote{But it is not in general transitive, compare \cite[Lemma~2.2]{confing}.} Two separations that are not nested are said to \emph{cross}.

A separation \AB\ is \emph{nested} with a set \cS\ of separations, written as ${\AB\| \cS}$, if $\AB \| \CD$ for every $\CD\in \cS$. A set \cS\ of separations is \emph{nested} with set~$\cS'$ of separations, written as $\cS\|\cS'$, if $\AB \| \cS'$ for every $\AB\in\cS$; then also $\CD\|\cS$ for every $\CD\in\cS'$.

A set of separations is called \emph{nested} if every two of its elements are nested. It is called \emph{symmetric} if whenever it contains a separation \AB\ it also contains \BA, and \emph{antisymmetric} if it contains no separation \AB\ together with its \emph{inverse}~\BA. A symmetric set of separations is a \emph{system} of separations, or \emph{\sys\/}; if all its separations are proper, it is a \emph{proper separation system}.\looseness=-1

A separation \AB\ \emph{separates} a set $X \sub V$ if $X$ meets both $A\sm B$ and $B\sm A$. Given a set \cS\ of separations, we say that $X$ is \emph{$\cS$-inseparable\/} if no separation in~$\cS$ separates~$X$. An \emph{\cS-block} of $G$ is a maximal $\cS$-inseparable set of vertices.%
   \COMMENT{}

Recall that a \emph{\td} of $G$ is a pair \TV\ of a tree $\cT$ and a family $\cV=(V_t)_{t\in \cT}$ of vertex sets $V_t\sub V(G)$, one for every node of~$\cT$, such that:
\begin{enumerate}[(T1)]
\item $V(G) = \bigcup_{t\in T}V_t$;
\item for every edge $e\in G$ there exists a $t \in T$ such that both ends of $e$ lie in~$V_t$;
\item $V_{t_1} \cap V_{t_3} \sub V_{t_2}$ whenever $t_2$ lies on the $t_1$--$t_3$ path in~$T$.
\end{enumerate}

The sets $V_t$ in such a \td\ are its \emph{parts}. Their intersections $V_t\cap V_{t'}$ for edges $tt'$ of the \emph{decomposition tree}~$\cT$ are the \emph{adhesion sets} of~\TV; their maximum size is the \emph{adhesion} of \TV.

Deleting an oriented edge $\vec e = t_1t_2$ of \cT\ divides $\cT-e$ into two components $T_1\owns t_1$ and ${T_2\owns t_2}$. Then $(\bigcup_{t\in T_1} V_t, \bigcup_{t\in T_2} V_t)$ is a separation of~$G$ with separator ${V_{t_1}\cap V_{t_2}}$ \cite[Lemma~12.3.1]{DiestelBook10noEE}; we say that our edge $\vec e$ \emph{induces\/} this separation. A node $t \in \cT$ is a \emph{hub node} if the corresponding part $V_t$ is the separator of a separation induced by an edge of \cT\ at $t$. If $t$ is a hub node, we call $V_t$ a \emph{hub}.

As is easy to check, the separations induced by (the edges of~\cT\ in) a \td~\TV\ are nested. Conversely, we proved in~\cite{confing} that every nested separation system%
   \footnote{In~\cite{confing}, as here, we needed this only for proper separation systems. However, the result and proof remain valid for arbitrary nested \sys s.}
    is induced by some \td:

\begin{thm}\label{treedec} {\rm \cite[Theorem~4.8]{confing}}
Every nested proper \sys~$\cN$ is induced by a \td\ $\left(\cT,\cV\right)$ of $G$ such that\vskip-3pt\vskip0pt
\begin{enumerate}[\rm (i)]\itemsep=0pt
\item every \cN-block of $G$ is a part of the decomposition;
\item every part of the decomposition is either an \cN-block of $G$ or a hub.
\end{enumerate}
\end{thm}

\noindent
See~\cite{confing} for how these \td s are constructed from~\cN.

\goodbreak

Let $k$ be a positive integer. A set $I$ of at least $k$ vertices is called \emph{\lk-insep\-ar\-able} if it is \cS-inseparable for the set $\cS=\cS_k$ of all separations of order~$<k$, that is, if for every separation $\AB\in\cS_k$ we have either $I \sub A$ or $I \sub B$. A~maximal \lk-inseparable set of vertices is a \emph{$k$-block} of~$G$.

Since a $k$-block is too large to be contained in the separator $A\cap B$ of a separation \AB\ of order~$<k$, it thus `chooses' one of the sides $A$ or~$B$, the one containing it. Equivalently, every $k$-block $b$ chooses one of each inverse pair \AB\ and \BA\ of separations in~$S_k$ if we want it to lie in the right-hand side of that separation. Let us give the set of separations chosen by~$b$ a name:
\begin{equation}\label{eq_Pb}
P_k(b) := \{\,\AB :  |A\cap B| < k\, \wedge\, b \sub B\,\}.
\end{equation}

Another way of making informed choices for small-order separations, but one that cannot necessarily be defined by setting a target set like~$b$, are \emph{tangles\/}. Tangles were introduced by Robertson and Seymour~\cite{GMX}; the definition can also be found in~\cite{DiestelBook10noEE}. Like $k$-blocks, $k$-tangles (those of order~$k$) have been considered as a way of identifying the `$k$-connected components' of a graph.

In order to treat blocks and tangles together in a unified way, let us distill the common essence of the `informed choices' they define for small-order separations into a couple of axioms, and then just work with these.%
   \COMMENT{}

One common property of $k$-tangles $P$ and the sets $P=P_k(b)$ defined by a $k$-block is
 $$\AB\in P\ \wedge\ \CD\le\AB\ \Rightarrow\ \DC\notin P. $$
   We shall call sets $P$ of separations that satisfy this implication \emph{consistent}.%
   \footnote{Note that we do not require that \CD\ must lie in~$P$, only that \DC\ does not. The term `consistent' is natural if we think of a separation \AB\ as `pointing towards'~$B$, so that all the elements of $P_k(b)$ point towards~$b$: then \DC\ points away from~$D$, which makes it inconsistent with~\AB, which points towards~$B$, if $\CD\le\AB$ and hence $D\supseteq B$.}
    Note that consistent sets of separations are antisymmetric: since $(A,B)\le (A,B)$, a~consistent set~$P$ cannot contain both $(A,B)$ and its inverse~$(B,A)$.%
   \COMMENT{}

Another common property of $k$-tangles $P$ and the sets $P=P_k(b)$ is this:
\medskip

 {\em \hfill For all $\AB,\CD \in P$ we have $(B\cap D,A\cup C) \notin P$.\hfill\rm (P)}%
   \COMMENT{}
   \medskip

If we think of the right-hand side $B$ of \AB, the side to which it `points', as its large side, condition~(P) becomes reminiscent of the property of (ultra)filters that `the intersection of large sets are large'. An important difference is that rather than demanding that $(A\cup C,B\cap D)\in P$, condition~(P) only%
   \COMMENT{}
   asks that the converse of this separation shall not be in~$P$.%
   \COMMENT{}
    Similarly, consistent sets of separations have a property reminiscent of an ultrafilter being non-principal:%
   \COMMENT{}
\begin{txteq}\label{ABandBA}
  \text{If $P$ is consistent, it contains no separation of the form~$(V,A)$}.
\end{txteq}
Indeed, as $(A,V)\le (V,A)$, a~consistent set of separations containing~$(V,A)$ must not contain the inverse of~$(A,V)$, which is $(V,A)$.%
   \footnote{If $P$ is an orientation of~$S_k$ (see below) and satifsies~(P), the condition that $(V,A)\notin P$ for all $A\sub V$ is in fact equivalent to the consistency of~$P$.}%
   \COMMENT{}
   Note that consistent sets of separations can contain improper separations of the form~$(A,V)$.

Note that while consistency says something only about nested separations, condition~(P) is essentially about crossing separations. A consistent set of separations satisfying~(P) will be called a \emph{profile}. Consistent sets of separations that do not satisfy~(P) will play a role later, too.%
   \footnote{As a typical example, consider the union of three large complete graphs $X_1,X_2,X_3$ identified in a common triangle. The three 3-separations whose left side is one of $X_1,X_2,X_3$ are consistent but do not satisfy~(P), because the separation $(B\cap D, A\cup C)$ in (P) happens to be one of the original three separations. The analogous system with four complete graphs does satisfy~(P).}%
   \COMMENT{}

Since consistency and condition~(P) only require the absence of certain separations from~$P$, a~requirement easy to meet by making $P$ small, profiles, unlike blocks and tangles, do not as such witness the existence of any highly connected substructure in a graph.%
   \footnote{Readers familiar with the notion of {\em preferences}, or {\em havens}~-- a way of making consistent choices of components of $G-X$ for small sets $X$ of vertices~-- will recognize this: it is because a preference or haven assigns a component of $G-X$ to {\em every} set $X$ of $<k$ vertices for some~$k$ that the bramble formed by these components has order~$\ge k$.}
    But they do as soon as we make them large, e.g.\ by requiring that they contain one of every pair of inverse separations in~$S_k$. To give such rich profiles a name, let us call a profile $P$ a \emph{$k$-profile} if it satisfies
\begin{txteq}\label{P_complete} 
Every separation in $P$ has order~$<k$, and for every separation \AB\ of order~$<k$ exactly one of \AB\ and \BA\ lies in~$P$.\looseness=-1
\end{txteq}%
   \COMMENT{}

By~\eqref{ABandBA}, every $k$-profile contains every separation $(A,V)$ with $|A|<k$.
As we have seen, the set $P_k(b)$ in~\eqref{eq_Pb} is a $k$-profile; we call it the $k$-profile \emph{induced by},  or simply {\em of},~$b$.%
   \COMMENT{}
   A $k$-profile induced by some $k$-block is a \emph{$k$-block profile}.

Since a $k$-block is a maximal \lk-inseparable set of vertices, there is for every pair of distinct $k$-blocks $b,b'$ a separation $\AB$ of order~$<k$ such that $\AB\in P_k(b)$ and $\BA \in P_k(b')$ \cite[Lemma~2.1]{confing}. Hence $P_k(b) \neq P_k(b')$. Thus, while every $k$-block induces a $k$-profile, conversely a $k$-profile $P$ is induced by at most one $k$-block, which we then denote by~$b(P)$. All $k$-block profiles $P$ then satisfy $P = P_k(b(P))$, and we say that $b$ and~$P$ \emph{correspond\/}.

Not every $k$-profile is induced by a $k$-block. For example, there are tangles of order~$k$ that are not induced by a $k$-block,%
   \COMMENT{}
   such as the unique tangle of any order $k\ge 5$ in a large grid%
   \COMMENT{}
   (which has no $k$-block for~$k\ge 5$; see~\cite[Example~3]{ForcingBlocks}). Conversely, there are $k$-block profiles that are not tangles; indeed, there are graphs that have interesting $k$-block profiles but have no non-trivial tangle at all \cite[Examples~4--5 and Section~6]{ForcingBlocks}. The notion of a $k$-profile thus unifies the ways in which $k$-blocks and tangles of order $k$ choose one side of every separation of order~$<k$, but neither of these two instances of $k$-profiles generalizes the other.

Let $\cS$ be any set of separations of~$G$. An \cS-block $X$ of $G$ is called \emph{large} (with respect to~\cS) if it is not contained in the separator of a separation in~$\cS$. If all the separations in $\cS$ have order~$< k$, an obvious but typical reason for an \cS-block to be large is that it has $k$ or more vertices. In analogy to \eqref{eq_Pb} we define for a large \cS-block~$X$
\begin{equation}\label{eq_PX}
P_{\cS}(X) := \{\AB \in \cS \mid X \sub B\}\sub\cS.
\end{equation}
Clearly, $P_{\cS}(X)$ is a profile; we call it the \emph{$\cS$-profile of $X$}. As before, the \cS-profiles $P_{\cS}(X)$ and $P_{\cS}(X')$ of distinct large \cS-blocks $X,X'$ are distinct.%
   \COMMENT{}

Not every $k$-profile has this form.%
   \COMMENT{}
   For example, a tangle $\theta$ of order ${k\ge 5}$ in a large grid is not the \cS-profile of a large \cS-block $X$ for any set $\cS\supseteq\theta$%
   \COMMENT{}
   of separations, since $X$ would be contained in%
   \COMMENT{}
   a large $\theta$-block but the grid has none.\looseness=-1 %
   \COMMENT{}

Although profiles are, formally, sets of separations, our intuition behind them is that they signify some `highly connected pieces' of our graph~$G$. Our aim will be to separate all these pieces in a tree-like way, and we shall therefore have to speak about sets of separations that, initially, are quite distinct from the profiles they are supposed to separate. To help readers keep their heads in this unavoidable confusion, we suggest that they think of the sets \cS\ of separations discussed below as (initially) quite independent of the profiles $P$ discussed along with them, the aim being to explore the relationship between the two.

A separation \AB\ \emph{distinguishes}%
   \COMMENT{}
   two subsets of $V$ if one lies in~$A$, the other in~$B$, and neither in~$A\cap B$. A set \cS\ of separations {\em distinguishes\/} two sets of vertices if some separation in~\cS\ does.

A separation \AB\ \emph{distinguishes} two sets $P,P'$ of separations if $\AB\in P\sm P'$ and $\BA\in P'\sm P$, or vice versa. Thus, a \lk-separation distinguishes two $k$-blocks if and only if it distinguishes their $k$-profiles. A set of separations~\cS\ \emph{distinguishes} $P$ from~$P'$ if some separation in \cS\ distinguishes them, and \cS\ \emph{distinguishes} a set \cP\ of sets of separations if it distinguishes every two elements of~\cP. By~\eqref{ABandBA},
 \begin{txteq}\label{distproper}
   Only proper separations can distinguish two consistent sets of separations, e.g., two profiles.
\end{txteq}%
   \COMMENT{}

An antisymmetric set $P$ of separations \emph{orients} a set \cS\ of separations if, for every $\AB\in\cS$, either ${\AB\in P}$ or $\BA\in P\cap\cS$. We then call $P\cap\cS$ an \emph{orientation\/} of~\cS, and an \emph{\cS-profile\/} if it%
   \COMMENT{}
   is a profile.

Distinct consistent orientations $O,O'$ of~$S$ are distinguished by some separation in~$S$.%
   \COMMENT{}
   Indeed, as $O$ and $O'$ are distinct we may assume that there is a separation $\AB\in O\sm O'$. Then $\BA\notin O$, since $O$ is consistent,%
   \COMMENT{}
  and $\BA\in O'$, since $\AB\in S$ and $O'$ orients~$S$. Hence \AB\ distinguishes $O$ from~$O'$.

If $P$ orients~$S$ and some set $X$ of vertices lies in $B$ for every $\AB\in P\cap\cS$, we say that $P$ \emph{orients} \cS\ \emph{towards\/}~$X$.%
   \COMMENT{}
   If $P$ is a profile then so is $P\cap\cS$; we call it the \emph{\cS-profile of~$P$\/}.%
   \footnote{This formalizes the idea that~$P$, thought of as a big chunk of~$G$, lies on exactly one side of every separation in~$\cS$. For example, if the separations in \cS\ have order~$<k$ and $\cS$ is symmetric, then every $k$-profile will orient~\cS.}

A profile orienting a set \cS\ of separations need not orient it towards any non-empty set of vertices: consider, for example, our earlier tangle~$\theta$ with $\cS = \theta$. However, a consistent set $P$%
   \COMMENT{}
   orienting a \emph{nested\/} set \cN\ of separations always orients it towards the intersection $X$ of all the sides~$B$ of separations $(A,B)\in {P\cap\cN}$.%
   \COMMENT{}
   By consistency, this set contains all the separators of the $\le$-maximal separations in~${P\cap\cN}$,%
   \COMMENT{}
   so it is non-empty if $G$ is connected.%
   \COMMENT{}
   If $X$ is an $\cN$-block, we say that $P$ \emph{lives in\/} this \cN-block~$X$.

However, $X$~need not be an \cN-block: although it is \cN-inseparable, it can be contained in a separator of a separation in~$\cN$%
   \COMMENT{}
   and extend to more than one \cN-block of~$G$;%
   \COMMENT{}
   then $P$ will not orient~\cN\ towards any of these. This can happen even if $P$ is a profile,%
   \footnote{For example, let $Z\sub V$ be such that $G-Z$ has four components $C_1,\dots,C_4$. Let $A_i = V(C_i)\cup Z$ and $B_i = V\sm V(C_i)$. Then $P = \{(A_i,B_i)\mid i=1,\dots,4\}$ is a nested profile, with $X = Z = A_i\cap B_i$ for all~$i$. }
   but not if $P$ is a $k$-profile: then~$X$ cannot be contained in the separator of any separation in~$P\cap\cN$.%
   \COMMENT{}
   This is interesting, because it makes $k$-profiles of nested separation systems `big' in a way arbitrary profiles need not be.%
   \COMMENT{}
   So if $P$ is a $k$-profile, then $X$ is an \cN-block towards which $P$ orients~\cN.%
   \COMMENT{}%
   \footnote{Conversely, if $P$ orients~\cN\ towards any \cN-block, then this is clearly equal to~$X$.}%
   \COMMENT{}

Given a set \cS\ of separations of $G$ and a set \cP\ of profiles orienting~\cS, let us say that two profiles $P,P'\in\cP$ \emph{agree on~\cS\/} if their \cS-profiles coincide, that is, if $P\cap\cS = P'\cap\cS$. This is an equivalence relation on~\cP, whose classes we call the \emph{\cS-blocks\/} of~\cP. By definition, elements $P,P'$ of the same \cS-block~\cQ\ of~\cP\ have the same \cS-profile $P\cap\cS = P'\cap\cS$, which we call the \emph{\cS-profile of~\cQ}.

A~separation~$\AB$ \emph{splits} a consistent set~$P$ of separations if both ${P\cup \{\AB\}}$ and $P \cup \{\BA\}$ are consistent. (This implies that neither \AB\ nor \BA\ is in~$P$.)%
   \COMMENT{}
   For example, the \cS-profile corresponding to an \cS-block $\cQ$ of a set \cP\ of profiles orienting a \sys~\cS\ is split by every separation \AB\ that distinguishes some distinct profiles in~$\cQ$.
By~\eqref{ABandBA}, every separation splitting a consistent set of separations must be proper.

We shall need the following lemma.

\begin{lem}Let \cN\ be a nested \sys.\label{splitlemma}\vskip-6pt\vskip-6pt
\begin{enumerate}[\rm (i)]\itemsep=-6pt
\item Every proper%
   \COMMENT{}
   separation $\AB \notin \cN$ that is nested with~\cN\ splits a unique consistent orientation~$O$ of~\cN.%
   \COMMENT{}
    This set $O$ is given by
  $$O = \{\CD \in \cN \mid \CD \le \AB\} \cup \{ \CD\in\cN \mid \CD\le\BA \}.$$
\item If two separations%
   \COMMENT{}
   not contained in but nested with~\cN\ split distinct consistent orientations of~\cN, they are nested with each other.
\end{enumerate}
\end{lem}

\begin{proof}
(i)~Since $\AB$ is nested with~\cN, for every separation $\CD \in \cN$ either $\CD$ or $\DC$ is smaller than one of $\AB$ or $\BA$ and thus contained in
$$O := \{\CD \in \cN \mid \CD \le \AB\} \cup \{ \CD\in\cN \mid \CD\le\BA \}.$$
By definition, $O$ contains only separations from~\cN. As we have seen,%
   \COMMENT{}
   every separation from \cN\ or its inverse lies in~$O$. Once we know that $O$ is consistent it will follow that $O$ is antisymmetric, so it will be an orientation of~\cN.

To check that $O$ is consistent, consider separations $\EF\le\CD$ with $\CD \in O$. Our aim is to show that $\FE\notin O$. This is clearly the case if $\EF\notin\cN$, since $O\sub\cN$ and $\cN$ is symmetric, so we assume that ${\EF\in\cN}$. By definition of~$O$, either $\CD\le\AB$ or $\CD\le\BA$; we assume the former. Then by transitivity $\EF\le\AB$, and hence $\EF \in O$ by definition of~$O$. To show that $\FE\notin O$ we need to check that $\FE\not\le\AB$ and $\FE\not\le\BA$. If $\FE\le\AB$ then $\BA\le\EF\le\AB$ and hence $B\sub A$, contradicting our assumption that \AB\ is proper. If $\FE\le\BA$ then $\AB\le\EF\le\AB$ and hence $\EF = \AB$, contradicting our assumption that $\AB\notin\cN$.

So $O$ is a consistent orientation of~\cN, and in particular antisymmetric. By the definition of~$O$, this implies that ${O \cup \{\AB\}}$ and $O \cup \{\BA\}$ are consistent.%
   \COMMENT{}
   Hence \AB\ splits~$O$, as desired.

It remains to show that $O$ is unique. Suppose \AB\ also splits a consistent orien\-ta\-tion $O' \neq O$ of~\cN. Let $\CD \in \cN$ distinguish $O$ from~$O'$, with $\CD\in O$ and ${\DC\in O'}$ say. By definition of~$O$, either $\CD\le\AB$ or ${\CD\le\BA}$. In the first case $O' \cup \{\AB\}$ is inconsistent, since ${\BA\le\DC\in O'\cup\{\AB\}}$ but also $\AB\in O'\cup\{\AB\}$. In the second case, $O'\cup\{\BA\}$ is inconsistent, since $\AB\le\DC\in O'\cup\{\BA\}$ but also $\BA\in O'\cup\{\BA\}$.

(ii)~Consider separations $\AB,(A',B') \notin \cN$ that are both nested with~\cN. Assume that \AB\ splits a consistent orientation $O$ of~\cN, and that $(A',B')$ splits a consistent orientation ${O' \neq O}$ of~\cN. From~\eqref{ABandBA} we know that \AB\ and $(A',B')$ must be proper separations, so they satisfy the premise of~(i) with respect to $O$ and~$O'$. As ${O \neq O'}$, there is a separation $\CD \in \cN$ with ${\CD \in O}$ and $\DC \in O'$. By the descriptions of $O$ and~$O'$ in~(i),%
   \COMMENT{}
   the separation \CD\ is smaller than \AB\ or~\BA, and \DC\ is smaller than $(A',B')$ or~$(B',A')$. The latter is equivalent to \CD\ being greater than $(B',A')$ or~$(A',B')$. Thus, $(B',A')$ or~$(A',B')$ is smaller than \CD\ and hence than \AB\ or~\BA, so \AB\ and $(A',B')$ are nested.
\end{proof}

We remark that the consistent set $O$ in Lemma~\ref{splitlemma}\,(i) is usually an \cN-profile; it is not hard to construct pathological cases in which $O$ fails to satisfy~(P), but such cases are rare.%
   \COMMENT{}

Finally, let us note that every \sys~\cS\ has a consistent orientation. Indeed, the consistent orientations of~\cS\ are precisely the orientations $P$ of~\cS\ that are closed down under~$\le$, those that contain $\CD\in\cS$ whenever $\CD\le\AB$ with $\AB\in P$. And we can obtain such an orientation of \cS\ greedily: consider any pair of proper separations $\{\AB,\BA\}\sub\cS\sm P$, and include in~$P$ both $\AB$ (say) and every separation $\CD\le\AB$ from~\cS. The resulting set $P$ is clearly closed down in~\cS\ under~$\le$, and remains so if we add all improper separations of the form~$(A,V)\in\cS$. It is easy to check that our collection of separations is now an orientation of~\cS, in particular, that we never included the inverse of a separation added to $P$ previously.%
   \COMMENT{}

\section{Tasks and strategies}\label{sec_Ops}

In this section we describe a systematic approach to distinguishing some or all of the $k$-profiles in~$G$ by (the separations induced by) canonical \td s of adhesion less than~$k$. Since the separations induced by a \td\ are nested, our main task in finding such a \td\ will be to select from the set $\cS$ of all proper%
   \COMMENT{}
   $(<k)$-separations of~$G$, which distinguishes the set of all $k$-profiles by~\eqref{ABandBA},%
   \COMMENT{}
   a nested subset \cN\ that will still distinguish all the $k$-profiles under consideration.

We begin by formalizing the notion of such `tasks'. We then show how to solve `feasible' tasks in various ways, and give examples showing how different strategies~-- all canonical in that they commute with graph iso\-mor\-phisms~-- can produce quite different solutions.

Consider a proper \sys~\cS\ and a set~\cP\ of profiles in~$G$. Let us call the pair $(\cS,\cP)$ a \emph{task\/} if every profile in \cP\ orients~\cS\ and \cS\ distinguishes~\cP. Another task $(\cS',\cP')$ is a \emph{subtask} of the task~$(\cS,\cP)$ if $\cS'\sub\cS$ and $\cP'\sub\cP$. 

The two conditions in the definition of a task are obvious minimum requirements which \cS\ and~\cP\ must satisfy before it makes sense to look for a nested subset $\cN\sub\cS$ that distinguishes~\cP. But to ensure that \cN\ exists, $\cS$~must also be rich enough (in terms of~\cP): the more profiles we wish to separate in a nested way, the more separations will we need to have available. For example, if \cS\ consists of two crossing separations \AB, \CD\ and their inverses, and $\cP$ contains the four possible orientations of~$\cS$ (which are clearly profiles), then \cS\  distinguishes~\cP\ but is not nested, while the two subsystems $\{\AB,\BA\}$ and $\{\CD,\DC\}$ of~\cS\ are nested but no longer distinguish~\cP. But if we enrich~\cS\ by adding two `corner separations' $(A\cap C, B\cup D)$, $(A\cup C, B\cap D)$ and their inverses, then these together with \AB\ and~\BA\ (say) form a nested subsystem that does distinguish~\cP.

More generally, we shall prove in this section that we shall be able to find the desired~$\cN$ if $\cS$ and \cP\ satisfy the following condition:
\begin{txteq}\label{eq_distwell}
Whenever $\AB,\CD \in \cS$ cross and there exist $P,P' \in \cP$ such that $\AB,\CD \in P$ and $\BA,\DC \in P'$, there exists a separation $\EF \in P\cap \cS$ such that $ (A\cup C,B\cap D) \le \EF$.
\end{txteq}
Anticipating our results, let us call a task $(\cS,\cP)$ \emph{feasible\/} if \cS\ and \cP\ satisfy~\eqref{eq_distwell}.

Let us take a moment to analyse condition~\eqref{eq_distwell}. Note first that, like the given separations \AB\ and~\CD, the new separation \EF\ will again distinguish $P$ from $P'$: by assumption we have $\EF \in P$, and by \eqref{eq_sym} we have ${\FE\le (B\cap D, A\cup C) \le \BA\in P'}$, so $\FE \in P'$ by the fact that $P'$ orients all of~\cS\ consistently.

Now the idea behind~\eqref{eq_distwell} is that in our search for~$\cN$ we may find ourselves facing a choice between two crossing separations $\AB,\CD\in\cS$ that both distinguish two profiles $P,P'\in\cP$, and wonder which of these we should pick for~\cN. (Clearly we cannot take both.) If~\eqref{eq_distwell} holds, we have the option to choose neither and pick $\EF$ instead: it will do the job of distinguishing $P$ from~$P'$, and since it is nested with both \AB\ and~\CD, putting it in~\cN\ entails no prejudice to any future inclusion of either \AB\ or~\CD\ in~\cN.

Separations in \cS\ that do not distinguish any profiles in~\cP\ are not really needed for~\cN, and so we may delete them.%
   \footnote{But do not have to: the freedom to discard or keep such separations will be our source of diversity for the \td s sought~-- which, as pointed out earlier, we may wish to endow with other desired properties than the minimum requirement of distinguishing~\cP.}%
   \COMMENT{}
   So let us call a separation \emph{\cP-relevant} if it distinguishes some pair of profiles in \cP, denote by \cR\ the set of all \cP-relevant separations in~\cS, and call $(\cR,\cP)$ the \emph{reduction\/} of~$(\cS,\cP)$. If $(\cS,\cP) = (\cR,\cP)$, we call this task \emph{reduced}. Since all the separations \AB, \CD, \EF\ in~\eqref{eq_distwell} are \cP-relevant, \cR~inherits~\eqref{eq_distwell} from~\cS\ (and vice versa):
\begin{txteq}
$(\cR,\cP)$ is feasible if and only if $(\cS,\cP)$ is feasible.
\end{txteq}

Consider a task~\SP,%
   \COMMENT{}
   fixed until Example~\ref{subtaskex} near the end of this section. Our aim is to construct~\cN\ inductively, adding  a few separations at each step. A~potential danger when choosing a new separation to add to~\cN\ is to pick one that crosses another separation that we might wish to include later. This can be avoided if we only ever add separations that are nested with all other separations in~\cS\ that we might still want to include in~\cN. So this will be our aim.\looseness=-1

At first glance, this strategy might seem both wasteful and unrealistic: why should there even be a separation in $\cS$ that we can choose at the start, one that is nested with all others? However, we cannot easily be more specific:%
   \COMMENT{}
   since we want our nested subsystem \cN\ to be canonical, we are not allowed to break ties between crossing separations without appealing to an invariant of $G$ as a criterion, and it would be hard to find such a criterion that applies to a large class of graphs without specifying this class in advance. But the strategy is also more realistic than it might seem. This is because the set of pairs of profiles we need to distinguish by separations still to be picked decreases as \cN\ grows. As a consequence, we shall need fewer separations in~\cS\ to distinguish them. We may therefore be able to delete from~\cS\ some separations that initially prevented the choice of a desired separation \AB\ for~\cN\ by crossing it, because they are no longer needed to distinguish profiles in what remains of~\cP, thus freeing \AB\ for inclusion in~\cN.

To get started, we thus have to look for separations \AB\ in $\cS$ that are nested with all other separations in~\cS. This will certainly be the case for \AB\ if, for every $\CD \in \cS$, we have either $\CD\le\AB$ or $\DC\le\AB$;%
   \footnote{This implies that \AB\ is maximal in~\cS, but only because we are assuming that all separations in~\cS\ are proper: improper separations \CD\ can satisfy $\DC < \AB < \CD$.\looseness=-1}%
   \COMMENT{}
   let us call such separations \AB\ \emph{extremal} in~\cS.

Consider distinct extremal separations \AB\ and \CD\ in~$S$. Since~$\le$ is a partial ordering, they cannot satisfy $\CD\le\AB\le\CD$. Hence, up to renaming, $\DC\le\AB$, which shows that $\{\AB,\CD\}$ is inconsistent:%
   \COMMENT{}
   \begin{txteq}\label{extremalinconsistent}
Distinct extremal separations are inconsistent. In particular, they cannot lie in the same profile.
\end{txteq}

If $\DC\le\AB$, as above, then $\AB\not\le\CD$, because $\DC\le\CD$ would imply that $D\sub C$ and hence that \CD\ is improper (contrary to our assumption about~$S$). But $\DC\le\AB$ also implies that $\CD\not\le\AB$, as otherwise $V=D\cup C\sub A$, and \AB\ would be improper. Hence,
\begin{txteq}
Distinct extremal separations are $\le$-incomparable.
\end{txteq}

Extremal separations always exist in a feasible task~\SP, as long as \cS\ contains no superfluous separations (which might cross useful ones):

\begin{lem}\label{lem_extremal} 
If \SP\ is feasible and reduced, then every $\le$-maxi\-mal element of~\cS\ is extremal in~\cS.
\end{lem}

\proof Let \AB\ be a maximal separation in~\cS, and let $\CD\in\cS$ be any other separation. If \AB\ is nested with \CD\ it is comparable with \CD\ or \DC. Hence either $\CD\le\AB$ or $\DC\le\AB$ by the maximality of~\AB, as desired. We may thus assume that \AB\ and \CD\ cross.

Since \SP\ is reduced, \AB\ and \CD\ each distinguish two profiles from~\cP. Pick $P\in\cP$ containing~\AB. Since $P$ orients all of~\cS, it also contains \CD\ or \DC; we assume it contains~\CD. Now pick $P'\in\cP$ containing~\DC. If also $\BA\in P'$, then by~\eqref{eq_distwell} there exists an $\EF\in\cS\cap P$ such that $\AB\le (A\cup C,B\cap D) \le \EF$. Since \AB\ and \CD\ cross, the first of these inequalities is strict, which contradicts the maximality of~\AB. Hence $\AB\in P'\cap P$. Since $\cP$ is reduced, there exists $P''\in\cP$ containing~\BA. Applying~\eqref{eq_distwell} to $P''$ and either $P$ or~$P'$,%
   \COMMENT{}
   we again find an $\EF > \AB$%
   \COMMENT{}
   that contradicts the maximality of~\AB.
   \endproof

Note that the proof of Lemma~\ref{lem_extremal} uses crucially that \SP\ is feasible.%
   \COMMENT{}

\begin{lem}\label{lem_P_AB} 
If \SP\ is reduced, then for every extremal separation $\AB$ in~\cS\ there is a unique profile $P_\AB \in \cP$ such that $\AB\in P_\AB$.
\end{lem}%

\proof As \SP\ is reduced, there is a profile $P\in\cP$ containing~\AB. Suppose there is another such profile $P'\in\cP$. Then $P$ and~$P'$ are distinguished by some $\CD\in\cS$, because \SP\ is a task.%
   \COMMENT{}
   Since \AB\ is extremal, we may assume that $\CD\le\AB$. Then the fact that \DC\ lies in one of $P,P'$ contradicts the consistency of this profile.\endproof

By Lemma~\ref{lem_P_AB}, \AB\ distinguishes $P_{\AB}$ from every other profile in~\cP.%
   \COMMENT{}

\medbreak

Let us call a profile $P$ orienting~\cS\ \emph{extremal\/} with respect to~\cS\ if it contains an extremal separation of~\cS.%
   \COMMENT{}
   This will be the greatest, and hence the only maximal, separation in~$P\cap\cS$.%
   \COMMENT{}

As we have seen, an extremal profile is distinguished from every other profile in~\cP\ by some separation~$\AB$ that is nested with all the other separations in~\cS; this makes \AB\ a good choice for~\cN. The fact that made \AB\ extremal, and hence nested with all other separations in~\cS, was its maximality in~\cS\ (Lemma~\ref{lem_extremal}). In the same way we may ask whether, given any profile $P\in\cP$ (not necessarily extremal), the separations that are $\le$-maximal in ${P\cap\cS}$ will be nested with every other separation in~\cS: these are the separations `closest to~$P$', much as $(A,B)$ was closest to~$P_\AB$ (although there can now be many such separations).

Let us prove that the following profiles have this property:

\begin{defi}
Call a profile $P$ orienting~\cS\ \emph{well separated} in~\cS\ if the set of $\le$-maximal separations in ${P\cap\cS}$ is nested.
\end{defi}

Note that extremal profiles are well separated.%
   \COMMENT{}

\goodbreak

\begin{lem}\label{wellseparated}
Given a profile $P$ orienting~\cS, the following are equivalent:\vskip-6pt\vskip0pt
\begin{enumerate}[\rm (i)]\itemsep=0pt
   \item $P$ is well separated in \cS.
   \item Every maximal separation in $P\cap \cS$ is nested with all of~\cS.
   \item For every two crossing separations $\AB,\CD \in P\cap\cS$ there exists a separation $\EF \in P\cap \cS$ such that $ (A\cup C,B\cap D) \le\EF$.
\end{enumerate}
\end{lem}

\begin{proof}
The implication (ii)$\to$(i) is trivial; we show (i)$\to$(iii)$\to$(ii).

(i)$\to$(iii):
Suppose that $P$ is well separated, and consider two crossing separations $\AB, \CD \in P\cap \cS$. Let $(A',B')\ge\AB$ be maximal in ${P\cap\cS}$. Suppose first that $(A',B')\|\CD$. This means that $(A',B')$ is $\le$-comparable with either \CD\ or~\DC. Since \AB\ is not nested with \CD\ we have $(A',B')\not\le\CD$ and $(A',B')\not\le\DC$, and since both $\CD$ and $(A',B')$ are in~$P$, its consistency yields $\DC\not\le(A',B')$. Hence $\CD\le(A',B')$, and thus $(A\cup C,B\cap D) \le (A',B')$.%
   \COMMENT{}
   This proves (iii) with $\EF:=(A',B')$.%
   \COMMENT{}

Suppose now that $(A',B')$ crosses~$\CD$. Let $(C',D')\ge\CD$ be maximal in $P\cap\cS$. Since $(A',B')$ and $(C',D')$ are both maximal in~$P\cap\cS$ they are nested, by assumption in~(i). As in the last paragraph, now with $(C',D')$ taking the role of $(A',B')$, and $(A',B')$ taking the role of~\CD, we can show that $\AB\le (A',B')\le (C',D')$ and hence  $(A\cup C,B\cap D) \le (C',D')$.%
   \footnote{In the argument we need that $(C,D)$ and $(A',B')$ cross. This is why we first treated the case that they don't (but in that case we used that \AB\ and \CD\ cross).}
   This proves (iii) with $\EF:=(C',D')$.%
   \COMMENT{}

(iii)$\to$(ii): Suppose some maximal $\AB$ in ${P\cap\cS}$ crosses some ${\CD \in \cS}$. As~$P$ orients~\cS, and by symmetry of nestedness, we may assume that ${\CD\in P}$. By~(iii), there is an $\EF \in P\cap \cS$ such that $ (A\cup C,B\cap D) \le \EF$, so $\AB\le\EF$ as well as $\CD\le\EF$. But then $\EF = \AB$ by the maximality of $\AB$, and hence $\AB\|\CD$, contradicting the choice of \AB\ and~\CD.
\end{proof}

Let us call a separation \AB\ \emph{locally maximal\/} in~\SP{}%
   \COMMENT{}
   if there exists a well-separated profile $P \in \cP$ such that \AB\ is $\le$-maximal in ${P\cap\cS}$. Lemma~\ref{wellseparated} shows that these separations are a good choice for inclusion in~\cN:

\begin{cor}\label{locmax}
Locally maximal separations in~\SP{}%
   \COMMENT{}
   are nested with all of~\cS.
\end{cor}

We have seen three ways of starting the construction of our desired nested subsystem $\cN\sub\cS$ by choosing for \cN\ some separations from \cS\ that are nested with all other separations in~\cS: we may choose either\vskip-3pt\vskip0pt

\begin{itemize}\itemsep=0pt%
   \COMMENT{}
\item the set $\sext\SP$ of extremal separations in~\cS\ and their inverses; or
\item the set $\sloc\SP$ of all locally maximal separations in~\SP\ and their inverses; or
\item the set $\smax\SP$ of all separations in~\cS\ that are nested with every separation in~\cS\ (which is a symmetric set).
\end{itemize}

\noindent
Clearly,
\begin{equation}\label{eq_contained}
\sext\SP \sub \smax\SP \supseteq \sloc\SP
\end{equation}
in general, and 
\begin{equation}\label{eq_contained_r}
\emptyset \neq \sext\SP \sub \sloc\SP \sub \smax\SP
\end{equation}
if $\cS\ne\es$ and $\SP$ is feasible and reduced,%
   \footnote{In fact, all we need for an extremal separation \AB\ to be locally maximal in a feasible \SP\ is that it lies in some $P\in\cP$. But this need not be the case if \SP\ is not reduced: although one of \AB\ and~\BA\ must lie in every $P\in\cP$ (because $P$ orients~\cS), it might happen that this is always~\BA.}
   since in that case every maximal separation in $\cS$ is extremal%
   \COMMENT{}
   (Lemma~\ref{lem_extremal}) and every extremal separation \AB\ is locally maximal for $P_\AB\in\cP$.%
   \COMMENT{}

\begin{ex}\label{ThreeBlobsEx}
Let $G$ consist of three large complete graphs $X_1,X_2,X_3$ threaded on a long path, as shown in Figure~\ref{ThreeBlobsFig}. Let \cS\ be the set of all proper 1-separa\-tions. Let $\cP = \{P_1,P_2,P_3\}$, where $P_i$ is the 2-profile induced by~$X_i$. Then $\smax\SP=\cS$, while \sloc\SP\ contains only the separations in \cS\ with separators $x_1, x_2, y_2$ and~$x_3$, and \sext\SP\ only those with separator $x_1$ or~$x_3$.
   \begin{figure}[htpb]
\centering
   	  \includegraphics{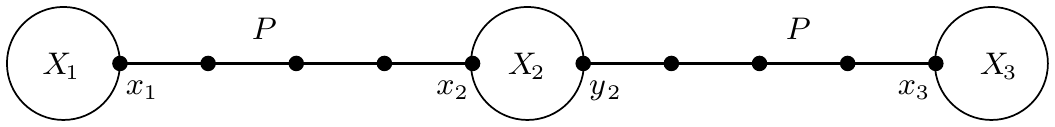}
   	  \caption{Different results for \sext\SP, \sloc\SP\ and~\smax\SP}
   \label{ThreeBlobsFig}\vskip-12pt\vskip0pt
   \end{figure}
\end{ex}

How shall we proceed now, having completed the first step of our algorithm by specifying some nested subsystem $\cN\in\{\sext\SP, \sloc\SP, \smax\SP\}$ of~\cS? The idea is that \cN~divides $G$ into chunks, which we now want to cut up further by adding more separations of~\cS\ to~\cN.%
   \COMMENT{}
   While it is tempting to think of those `chunks' as the \cN-blocks of~$G$, it turned out that this fails to capture some of the more subtle scenarios. Here is an example:

   \begin{figure}[htpb]
\centering
   	  \includegraphics{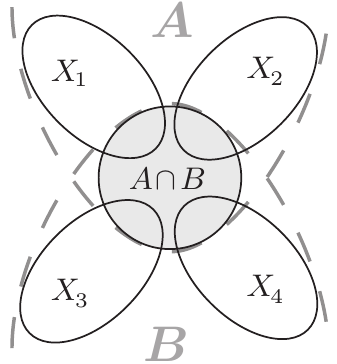}
   	  \caption{Two \cS-distinguishable profiles living in an \cS-inseparable \cN-block}
   \label{UnsplittableChunk}\vskip-12pt\vskip0pt
   \end{figure}

\begin{ex}\label{ABchunk}
Let $G$ be the graph of Figure~\ref{UnsplittableChunk}. Let $\cN$ consist of the separations $(X_1,Y_1),\dots,(X_4,Y_4)$ and their inverses $(Y_i,X_i)$, where $Y_i := (A\cap B)\cup\bigcup_{j\ne i} X_j$, and let ${\cS := \cN\,\cup\{\AB,\BA\}}$. Let $\cP$ consist of the following six profiles: the orientations of \cS\ towards $X_1,\dots,X_4$, respectively, and two further profiles $P$ and $P'$ which both orient \cN\ towards $A\cap B$ but of which $P$ contains~\AB\ while $P'$ contains~\BA. Then \cN\ distinguishes all these profiles except $P$ and~$P'$. But these are distinguished by~\AB\ and~\BA, so we wish to add these separations to~\cN.

The profiles $P$ and~$P'$ live in the same \cN-block of~$G$, the set $A\cap B$. But although \cS\ distinguishes $P$ from~$P'$, it does not separate this \cN-block. We therefore cannot extend~\cN\ to a \sys\ distinguishing~\cP\ by adding only separations from~\cS\ that separate an \cN-block of~$G$.
\end{ex}

The lesson to be learnt from Example~\ref{ABchunk} is that the `chunks' into which \cN\ divides our graph $G$ should not be thought of as the \cN-blocks of~$G$. An alternative that the example suggests would be to think of them as the \cN-blocks of~\cP: the equivalence classes of \cP\ defined by how its profiles orient~\cN.
   In the example, \cP~has five \cN-blocks: the four singleton \cN-blocks consisting just of the profile $P_i$ that orients \cN\ towards~$X_i$, and another \cN-block $\cQ = \{P,P'\}$. So the algorithm could now focus on the subtask $(\cR_\cQ,\cQ)$ with $\cR_\cQ = \{\AB,\BA\}$ consisting of those separations from~\cS\ that distinguish profiles in~\cQ.

More generally, we could continue our algorithm after finding \cN\ by iterating it with the subtasks $(\cR_\cQ,\cQ)$ of~\SP, where \cQ\ runs over the non-trivial \cN-blocks of~\cP\ and $\cR_\cQ$ is the set of \cQ-relevant separations in~\cS. This would indeed result in an overall algorithm that eventually produces a nested subsystem of~\cS\ that distinguishes~\cP, solving our task~\SP.

However, when we considered our three alternative ways of obtaining~\cN, we also had a secondary aim in mind: rather than working with the reduction \RP\ of~\SP\ straight away,\vadjust{\penalty-200} we kept our options open to include more separations in \cN\ than distinguishing~\cP\ requires, in order perhaps to produce a \td\ into smaller parts.%
   \footnote{In Example~\ref{ThreeBlobsEx} with~\sext\SP, where \cN\ consists of the proper 1-separations with separator $x_1$ or~$x_3$, every $\cN$-block of~\cP\ is trivial. But the middle \cN-block of $G$ consists of $X_2$ and the entire path~$P$, so we might cut it up further using the remaining 1-separations in~\cS. If \cP\ consisted only of $P_1$ and~$P_3$, then \sext\SP\ would have produced the same~\cN, and the middle \cN-block would not even have a profile from~\cP\ living in it. But still, we might want to cut it up further.\looseness=-1}
   In the same spirit, our secondary aim now as we look for ways to continue our algorithm from~\cN\ is not to exclude any separation of $\cS\sm\cN$ from possible inclusion into~\cN\ without need, i.e., to subdivide \SP\ into subtasks $(\cS_i,\cP_i)$ if possible with $\bigcup_i \cS_i = \cS$.%
   \COMMENT{}

In view of these two aims, the best way to think of the chunks left by~\cN\ turned out to be neither as the (large)%
   \COMMENT{}
   \cN-blocks of~$G$, nor as the \cN-blocks of~\cP, but as something between the two: as the set $\cO_\cN$ of all consistent orientations of~\cN. Let us look at these in more detail.

Recall that since every $P\in\cP$ orients~\cN, it defines an \cN-profile~$P\cap\cN$. Equi\-valent $P,P'$ define the same \cN-profile $P\cap\cN = P'\cap\cN$, the \cN-profile of the \cN-block \cQ\ containing them. This is a consistent orientation of~\cN. Conversely, given $O\in\cO_\cN$, let us write $\cP_O$ for the set of profiles $P\in\cP$ with $P\cap\cN = O$. Note that $\cO_\cN$ may also contain consistent orientations~$O$ of~\cN, including \cN-profiles, that are not induced by any $P\in\cP$, i.e., for which $\cP_O = \es$ (Fig.~\ref{ThreeBlobsFig}).

Similarly, every large \cN-block~$X$ of~$G$ defines an \cN-profile, the \cN-profile $P_{\cN\,} (X)$ of~$X$. This is a consistent orientation of~\cN. Again, $\cO_\cN$ may also contain consistent orientations that are not of this form.%
   \footnote{In Example~\ref{ABchunk}, the set $A\cap B$ is a small \cS-block of~$G$ for the nested \sys~\cS. The profiles $P,P'$ are two consistent orientations of \cS\ orienting it towards~$A\cap B$, but not towards any large \cS-block.}

Recall that a separation \AB\ \emph{splits} $O\in\cO_\cN$ if both $O\cup\{\AB\}$ and $O\cup\{\BA\}$ are again consistent.%
   \footnote{In Example~\ref{ABchunk}, the \cN-profile of $X=A\cap B$ could be split into the consistent sets $P$ and~$P'$ by adding the separations \AB\ and~\BA, although the large \cN-block $X$ could not be separated by any separation in~\cS. Thus, splitting the \cN-profile of a large \cN-block is more subtle than separating the \cN-block itself.
   \endgraf Although all the consistent sets of separations considered in this example are in fact profiles, our aim to retain all the separations from $\cS\sm\cN$ at this state requires that we do not restrict $\cO_\cN$ to profiles: there may be separations in~$\cS$ (which we want to keep) that only split a consistent orientation of \cN\ that is not a profile, or separations that split an \cN-profile into two consistent separations that are not profiles.}
    Let us write $\cS_O$ for the set of separations in~\cS\ that split~$O$. These sets $\cS_O$ extend our earlier sets $\cR_\cQ$ in a way that encompasses all of~$\cS\sm\cN$, as intended:

\begin{lem}\label{focus}
Let \cN\ be a nested \sys\ that is oriented by every profile in~\cP\ and is nested with~\cS.%
   \footnote{For better applicability of the lemma later, we do not require that $\cN\sub\cS$.}
\begin{enumerate}[\rm (i)]\itemsep=0pt\vskip-3pt\vskip0pt
   \item $(\,\cS_O\mid O\in\cO_\cN\,)$ is a partition of~$\cS\sm\cN$ (with $\cS_O=\es$ allowed).
   \item ${(\,\cP_O\mid O\in\cO_\cN\,)}$ is a partition of~\cP\ (with $\cP_O=\es$ allowed). 
   \item The \cN-profile~$P$ of any \cN-block \cQ\ of~\cP\ satisfies $\cP_P = \cQ$ and $\cS_P\supseteq\cR_\cQ$.
   \item If \SP\ is feasible, then the $(\cS_O,\cP_O)$ are feasible tasks.
\end{enumerate}
\end{lem}

\begin{proof}
(i) By Lemma~\ref{splitlemma}, every separation $\AB\in\cS\sm\cN$ splits a unique consistent orientation of~\cN. Note that \AB\ is proper because \SP\ is a task.

(ii) follows from the fact that every profile in~\cP\ orients~\cN\ consistently.

(iii) The first assertion is immediate from the definition of an \cN-block of~\cP. For the second assertion let $\AB\in\cR_\cQ$ be given, distinguishing $Q,Q'\in\cQ$ say. By~(i), we have $\AB\in\cS_O$ for some $O\in\cO_\cN\,$. Since $Q$ and~$Q'$ are consistent, agree with $P$ on~\cN, and orient $\{\AB,\BA\}$ differently, \AB~splits the consistent orientation $P$ of~\cN. By the uniqueness of $O$%
   \COMMENT{}
   this implies~$P=O$. Hence, $\AB\in\cS_O=\cS_P$ as desired.

(iv) As $\cS_O$ distinguishes~$\cP_O$, by~(iii), we only have to show that $(\cS_O,\cP_O)$ is feasible.%
   \COMMENT{}
   As \SP\ is feasible, there is a separation \EF\ in~\cS\ for any two crossing separations $\AB,\CD\in \cS_O$ distinguishing profiles $P,P'\in\cP_O$ as in~\eqref{eq_distwell}. Since \EF\ also distinguishes $P$ from~$P'$,%
   \COMMENT{}
   we have $\EF\in\cS_O$ by~(iii).\looseness=-1
\end{proof}

\noindent
We remark that the inclusion in Lemma~\ref{focus}\,(iii) can be strict, since $\cS_O$ may contain separations that do not distinguish any profiles in~\cP. Similarly, we can have $\cS_O\ne\es$ for $O\in\cO_\cN$ with $\cP_O=\es$.

\medskip
The subtasks $(\cS_O,\cP_O)$ will be `easier' than the original task~\SP, because we can reduce them further:

\begin{ex}\label{subtaskex}
The separations $\XY$ and $(X',Y')$ in Figure~\ref{subtaskfig} are \cP-relevant (because they separate the profiles $P,P'\in\cP$, say), so they will not be deleted when we reduce~\cS\ (which is, in fact reduced already). They both belong to~$\cS_O$ for the middle consistent orientation $O$ of~\cN, but are no longer $\cP_O$-relevant, where $\cP_O = \{P_1,P_2,P_3\}$ as shown. We can therefore discard them when we reduce the subtask $(\cS_O,\cP_O)$ before reapplying the algorithm to it, freeing \AB\ and~\CD\ for adoption into~\cN\ in the second step.
   \begin{figure}[htpb]
\centering
   	  \includegraphics{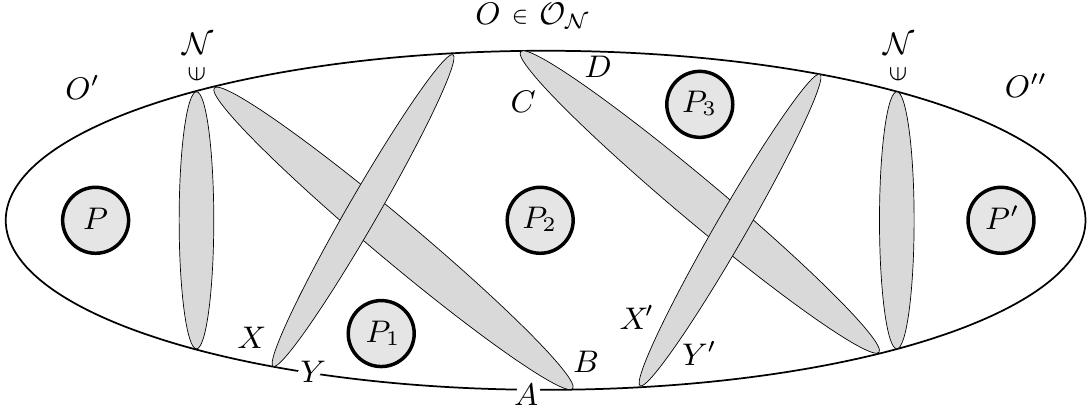}
   	  \caption{\XY\ and $(X',Y')$ are \cP-relevant but no longer $\cP_O$-relevant}
   \label{subtaskfig}\vskip-12pt\vskip0pt
   \end{figure}
\end{ex}

More generally, reducing a subtask $(\cS',\cP')$ will be the crucial step in getting our algorithm back afloat if it finds no separation in $\cS'$ that is nested with all the others. Example~\ref{subtaskex} shows that this can indeed happen.%
   \footnote{More generally, if we apply \smax\SP\ in the first step to obtain~\cN, say, then every subtask $(\cS_O,\cP_O)$ with $O\in\cN_O$ will have this property: if a separation $\AB\in\cS_O$ was nested with all of~$\cS_O$ it would in fact be nested with all of~$\cS$ (and have been included in~\cN\,), by Lemma~\ref{splitlemma}\,(ii).}
   But after reducing $(\cS',\cP')$ to $(\cR',\cP')$, say, we know from \eqref{eq_contained_r} that each of \nops\ will find a separation in~$\cR'$ that is nested with all the others.

As notation for the double step of first reducing a task ${\SP}$, to~${\RP}$ say, and then applying \sext, \sloc\ or~\smax, let us define%
   \footnote{For the remainder of this section, $G$ and \SP\ will no longer be fixed.}
  $$\srext{\SP} : = \sext{\RP};\quad \srloc{\SP} : = \sloc{\RP};\quad \srmax{\SP} : = \smax{\RP}.$$

\noindent
We shall view each of \ops\ as a function that maps a given graph $G$ and a feasible task ${\SP}$ in $G$ to a nested subsystem $\cN\,'\!$ of~$\cS'$.

\sloppy
A \emph{strategy}%
   \COMMENT{}
   is a map $\sigma\colon \N \to \{\ops\}$ such that ${\sigma(i) \in \{\rops\}}$ for infinitely many~$i$. The idea is that, starting from some feasible task~\SP, we apply $\sigma(i)$ at the $i$th step of the algorithm to the subtasks produced by the previous step, adding more and more separations to~\cN. The requirement that for infinitely many~$i$ we have to reduce the subtasks first ensures that we cannot get stuck before \cN\ distinguishes all of~$\cP$. 

\fussy
Formally, we define a map $(\sigma, G, \SP)\mapsto \cN_\sigma\SP$ by which every strategy~$\sigma$ \emph{determines} for every feasible task \SP\ in a graph~$G$ some set $\cN_\sigma\SP$. We define this map recursively, as follows. Define $\sigma^+$ by setting $\sigma^+(i) := \sigma(i+1)$ for all $i \in \N$. Note that if $\sigma$ is a strategy then so is~$\sigma^+$.%
   \COMMENT{}
   Let $s:=|\cS|$, and let $r_\sigma$ be the least integer~$r$ such that $\sigma(r) \in \{\rops\}$. Our recursion is on~$s$, and for fixed~$s$ on~$r_\sigma$, for all~$G$.

If $s=0$, we let $\cN_\sigma\SP:=\cS=\es$. Suppose now that $s\ge 1$; thus, $\cS\ne\es$. Let $\cN := \sigma(0)\SP$. By Lemma~\ref{focus}\,(iv), the subtasks $(\cS_O,\cP_O)$ with $O\in\cO_\cN$ are again feasible.

Assume first that $r_\sigma = 0$, i.e.\ that $\sigma(0) \in \{\rops\}$, and let \RP\ be the reduction of~\SP. If $\cR\subsetneq\cS$ we let $\cN_\sigma\SP := \cN_\sigma\RP$, which is already defined.%
   \COMMENT{}
   If $\cR = \cS$ then $\cN\ne\es$ by~\eqref{eq_contained_r}, and $|\cS_O|\le |\cS\sm\cN| < s$ for every $O\in\cO_\cN\,$. Thus, $\cN_{\sigma^+} (\cS_O,\cP_O)$ is already defined.

Assume now that $r_\sigma > 0$, i.e.\ that $\sigma(0) \in \{\nops\}$. Then $r_{\sigma^+} < r_\sigma$, so again $\cN_{\sigma^+} (\cS_O,\cP_O)$ is already defined. In either case we let
\begin{equation}\label{Nsigma}
\cN_\sigma\SP :=\ \cN\>\cup\!\bigcup_{O\,\in\,\cO_\cN}\!\! \cN_{\sigma^+} (\cS_O,\cP_O)\,.
\end{equation}

\begin{thm}\label{thm_can_sol}
Every strategy $\sigma$ determines for every feasible task~\SP\ in a graph~$G$ a~nested subsystem $\cN_\sigma$ of~$\cS$ that distinguishes all the profiles in~\cP.%
   \COMMENT{}

These sets $\cN_\sigma$ are canonical in that, for each~$\sigma$, the map $(G,\cS,\cP) \mapsto \cN_\sigma$ commutes with all isomorphisms $G\mapsto G'$.%
   \COMMENT{}
   In particular, if $\cS$ and $\cP$ are invariant under the automorphisms of~$G$, then so is~$\cN_\sigma$.%
   \COMMENT{}
\end{thm}

\proof
We apply induction along the recursion in the definition of $\cN_\sigma = \cN_\sigma\SP$. If $s=0$, then $\cN_\sigma=\cS$ distinguishes all the profiles in~\cP, because \SP\ is a task.\looseness=-1

Suppose now that $s\ge 1$. Then $\cN_\sigma$ is defined by~\eqref{Nsigma}. Both $\cN$ and the sets $\cN_{\sigma^+}(\cS_O, \cP_O)$ are subsets of~\cS, hence so is~$\cN_\sigma$. By definition, \cN\ is nested with all of~\cS, in particular, with itself and the sets $\cN_{\sigma^+}(\cS_O, \cP_O)$. These sets are themselves nested by induction, and nested with each other by Lemma~\ref{splitlemma}\,(ii). Thus, $\cN_\sigma$~is a nested subset of~\cS.

Any two profiles in the same \cN-block of~\cP\ are, by induction, distinguished by $\cN_{\sigma^+}(\cS_O, \cP_O)$ for their common consistent orientation~$O$ (cf.\ Lemma~\ref{focus}\,(iii)). Profiles from different \cN-blocks of \cP\ are distinguished by~\cN. Hence $\cN_\sigma$ distinguishes~\cP.

Finally, the maps $\SP\mapsto \cN_\sigma$ commute with all isomorphisms $G\mapsto G'$. Indeed, the maps $\SP\mapsto\cN$ and hence ${\SP\mapsto \{\,(\cS_O,\cP_O)\mid O\in\cO_\cN\,\}}$%
   \COMMENT{}
   do by definition of~\ops, and the maps $(\cS_O,\cP_O)\mapsto\cN_{\sigma^+}(\cS_O,\cP_O)$ do by induction.
\endproof

Let us complete this section with an example of how the use of different strategies can yield different nested \sys s. Unlike in the simpler Example~\ref{ThreeBlobsEx}, these will not extend each other, but will be incomparable under set inclusion. Let $\sExt$, $\sLoc$ and $\sMax$ denote the strategies given by setting $\sExt(i) = \srext$ and $\sLoc(i) = \srloc$ and $\sMax(i) = \srmax$, respectively, for all~$i\in\N$.

\begin{figure}[h]
  \begin{subfigure}{\textwidth}
    \includegraphics[width=\textwidth]{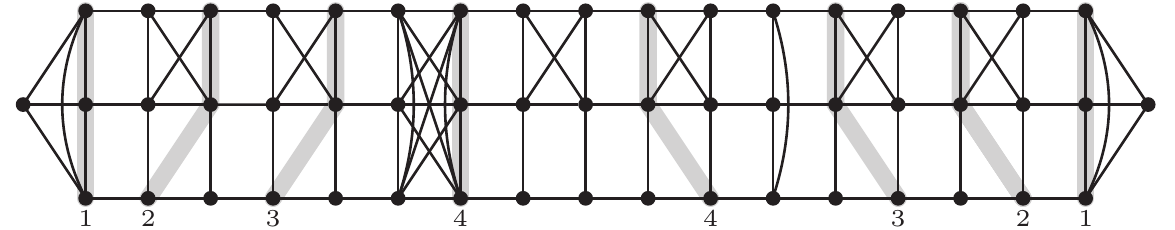}\vskip-3pt
    \caption*{\sExt: add the extremal separations at each step}
  \end{subfigure}\endgraf\bigskip
  \begin{subfigure}{\textwidth}
    \includegraphics[width=\textwidth]{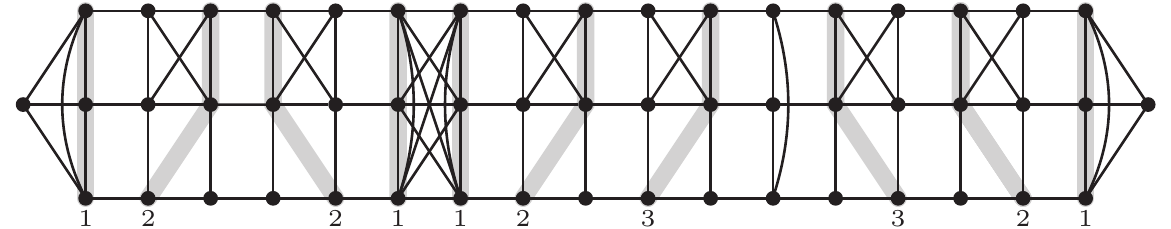}\vskip-3pt
    \caption*{\sLoc: add the locally maximal separations at each step}
  \end{subfigure}\endgraf\bigskip
  \begin{subfigure}{\textwidth}
    \includegraphics[width=\textwidth]{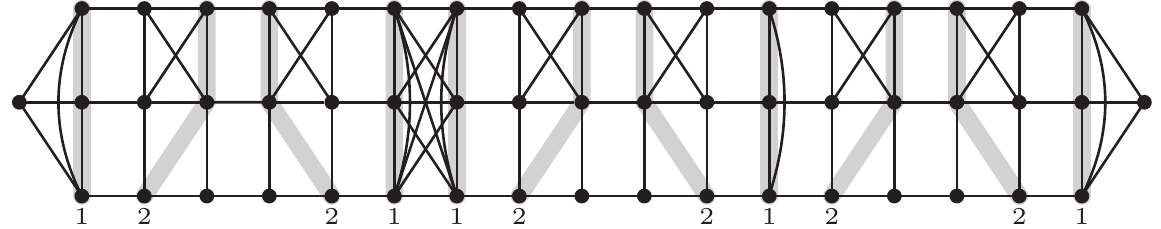}\vskip-3pt
    \caption*{\sMax: add all possible separations at each step}
  \end{subfigure}
  \caption{Three different nested separation systems distinguishing the $4$-blocks}\label{fig_theEx}\vskip-12pt\vskip0pt
\end{figure}

\begin{ex}\label{the_Example}
Let $G$ be the $3$-connected graph obtained from a $(3\times 17)$-grid by attaching two $K^4$s at its short ends, and some further edges as in Figure~\ref{fig_theEx}. Let \cS\ be the set of all its proper $3$-separations, and \cP\ the set of all its $4$-block profiles. It is not hard to show, and will follow from Lemma~\ref{lem_stdprob}, that \SP\ is a feasible task.

The grey bars in each of the three copies of the graph highlight the separators of the separations in $\cN_{\sExt} \SP$, in~$\cN_{\sLoc} \SP$, and in~$\cN_{\sMax} \SP$. The step at which a separator was added is indicated by a number.

Note that the three nested separation systems obtained are not only $\sub$-in\-comp\-arable. They are not even nested with each other: for every pair of $\cN_\sExt$, $\cN_\sLoc$ and~$\cN_\sMax$ we can find a pair of crossing separations, one from either system.
\end{ex}

\section{\boldmath Iterated strategies and \td s}\label{sec_iterate}

Let us apply the results of Section~\ref{sec_Ops} to our original problem of finding, for any set \cP\ of $k$-profiles in~$G$, within the set \cS\ of all proper $(<k)$-separations of $G$ a nested subset that distinguishes~\cP\ (and hence gives rise to a \td\ of adhesion $<k$ that does the same). This is easy now if $G$ is $(k-1)$-connected:

\begin{lem}\label{lem_stdprob}
If $G$ is $(k-1)$-connected $(k\ge 1)$, then \SP\ is a feasible task.
\end{lem}

\begin{proof}
By \eqref{ABandBA}, the pair \SP\ is a task: since all $k$-profiles contain the same improper separations, the separations $(A,V)$ with $|A|<k$, distinct $k$-profiles must differ on~$S$ and are thus distinguished by~$S$. So let us show that this task is feasible~\eqref{eq_distwell}.
Let \AB\ and \CD\ be crossing separations in~$\cS$, and let $P, P' \in \cP$ be such that $\AB,\CD \in P$ and $\BA,\DC \in P'$. We first prove that $\EF := (A\cup C, B \cap D)$ has order at most $k-1$.

Suppose $\EF$ has order greater than $k-1$. By Lemma~\ref{counting} this implies that the separation $\XY:= (B \cup D, A \cap C)$ has order less than $k-1$, and hence is improper since $G$ is $(k-1)$-connected. As both \AB\ and \CD\ are proper separations and hence $X\sm Y = (B\cup D)\sm (A\cap C)\ne\es$, we then have $Y\sub X$. Then $\XY\notin P'$, by \eqref{ABandBA}.%
   \COMMENT{}
   But by definition of \XY\ and (P) for~$P'$ we also have ${\YX \notin P'}$. This contradicts~\eqref{P_complete}.

We have shown that \EF \ has order at most $k-1$. By \eqref{P_complete} and~(P) this implies that $\EF\in P$.%
   \COMMENT{}
   To complete our proof of~\eqref{eq_distwell}, it remains to show that $\EF\in\cS$, i.e.\ that \EF\ is proper. If it is improper%
   \COMMENT{}
    then $F\sub E = V$, since $E\sm F\supseteq A\sm B\ne\es$ and therefore $F\ne V$. By~\eqref{ABandBA}, this contradicts the fact that $\EF\in P$.
\end{proof}

Coupled with Lemma~\ref{lem_stdprob}, we can apply Theorem~\ref{thm_can_sol} as follows:

\begin{cor}\label{cor_kcon}
Every strategy $\sigma$ determines for every $(k-1)$-connected graph~$G$ a canonical nested system of separations of order~$k-1$ which distinguishes all the $k$-profiles in~$G$. \qed
\end{cor}

If $G$ is not $(k-1)$-connected, the task%
   \COMMENT{}
   \SP\ consisting of the set \cS\ of all proper $(<k)$-separations of $G$ and the given set \cP\ of $k$-profiles need not be feasible. Indeed, the separation $(A\cup C,B\cap D)$ in~\eqref{eq_distwell} might have order $\ge k$ even if both \AB\ and \CD\ have order~$<k$. Then if $B\cap D$ induces a big complete graph, for example, there will be no \EF\ as required in~\eqref{eq_distwell}.

However, if $|A\cap B| = |C\cap D| = k-1$ in this example, the separation ${(B\cup D, A\cap C)}$ will have some order $\ell < k-1$. This separation, too, distinguishes the profiles $P,P'$ given in~\eqref{eq_distwell}, and thus by~\eqref{distproper} lies in~$S$. Hence our dilemma of having to choose between \AB\ and~\CD\ for inclusion in our nested subset \cN\ of~\cS\ (which gave rise to~\eqref{eq_distwell} and the notion of feasibility) would not occur if we considered lower-order separations first: we would then have included ${(B\cup D, A\cap C)}$ in \cN, and would need neither \AB\ nor \CD\ to distinguish $P$ from~$P'$.

It turns out that this approach does indeed work in general. Given our set~$\cP$ of $k$-profiles, let us define for any $1\le\ell\le k$ and $P\in\cP$ the induced $\ell$-profile
 $$P_\ell := \{\,\AB\in P : |A\cap B| < \ell\,\},\quad\text{and set}\quad \cP_\ell := \{\, P_{\ell} \mid P \in \cP\, \}\,.$$
Note that distinct $k$-profiles $P$ may induce the same $\ell$-profile~$P_\ell$. Let $\kappa(P,P')$ denote the least order of any separation in~$G$ that distinguishes two profiles $P,P'$.%
   \COMMENT{}

The idea now is to start with a nested set $\cN_1\sub\cS$ of $(<1)$-separations that distinguishes~$\cP_1$,%
   \COMMENT{}
    then to extend $\cN_1$ to a set $\cN_2\sub\cS$ of $(<2)$-separations that distinguishes~$\cP_2$,%
   \COMMENT{}
   and so on. The tasks $(\cS_O,\cP_O)$ to be solved at step~$k$, those left by the consistent orientations $O$ of $\cN_{k-1}$, will then be feasible: since $\cN_{k-1}$ already distinguishes~$\cP_{k-1}$, and hence distinguishes any $P,P'\in\cP$ with $\kappa(P,P') < k-1$,  any~$P,P'$ in a common~$\cP_O$ will satisfy $\kappa(P,P') = k-1$, and \eqref{eq_distwell} will follow from Lemma~\ref{counting} as in the proof of Lemma~\ref{lem_stdprob}.%
   \COMMENT{}

What is harder to show is that those $(\cS_O,\cP_O)$ are indeed tasks: that $\cS_O$ is rich enough to distinguish~$\cP_O$.%
   \COMMENT{}
   This will be our next lemma. Let us say that a separation \AB\ of order $\ell$ that distinguishes two profiles $P$ and~$P'$ does so \emph{efficiently} if $\kappa(P,P') = \ell$. We say that \AB\ is \emph{\cP-essential} if it efficiently distinguishes some pair of profiles in~\cP. Note that for ${\ell \le m}$ we have $(\cP_{m})_{\ell} = \cP_\ell$, and if \AB\ is $\cP_\ell$-essential it is also $\cP_m$-essential.

\begin{lem}\label{lem_sopo_feasible}
Let \cP\ be a set of $k$-profiles in~$G$, let \cN\ be a nested system of $\cP_{k-1}$-essential separations of~$G$ that distinguishes all the profiles in~$\cP_{k-1}$ efficiently, and let $\cS$ be the set of all proper $(k-1)$-separations of $G$ that are nested with~\cN. Then for every consistent orientation~$O$ of~\cN\ the pair \SPi O\ is a feasible task.
\end{lem}

\begin{proof}
As pointed out earlier,%
   \COMMENT{}
   \SPi O\ will clearly be feasible once we know it is a task. Since all profiles in~$\cP_O$ are $k$-profiles and hence orient~$\cS_O$, we only have to show that $\cS_O$ distinguishes~$\cP_O$. 

So consider distinct profiles $P_1,P_2 \in \cP_O$. Being $k$-profiles, they are distinguished by a proper%
   \COMMENT{}
   separation~\AB\ of order at most~$k-1$. Choose \AB\ nested with as many separations in~\cN\ as possible; we shall prove that it is nested with all of~\cN, giving $\AB\in\cS_O$ as desired. Note that $|A\cap B| = k-1$, since \cN\ does not distinguish  ${P_1,P_2\in\cP_O}$. As \AB\ distinguishes $P_1$ from~$P_2$, we may assume $\BA \in P_1$ and $\AB \in P_2$.

Suppose \AB\ crosses a separation $\CD \in \cN$. Since every separation in~\cN\ is $\cP_{k-1}$-essential, by assumption, there are profiles $Q'_1, Q'_2 \in \cP_{k-1}$ such that \CD\ distinguishes $Q'_1$ from $Q'_2$ efficiently. By definition of $\cP_{k-1}$, this implies that there are distinct profiles $Q_1, Q_2 \in \cP$ which \CD\ distinguishes efficiently. Then
\begin{equation}\label{ABorderCD}
   m:= |C\cap D| < k-1 = |A\cap B|.
\end{equation}
Hence \CD\ does not distinguish $P_1$ from~$P_2$;%
   \COMMENT{}
   we assume that $\DC \in P_1 \cap P_2 \cap Q_1$ and $\CD \in Q_2$. Since $Q_2$ is a $k$-profile it contains precisely one of $\AB$ and $\BA$, we assume $\AB \in Q_2$ (Figure~\ref{OptNestFig}).

   \begin{figure}[htpb]
\centering
   	  \includegraphics{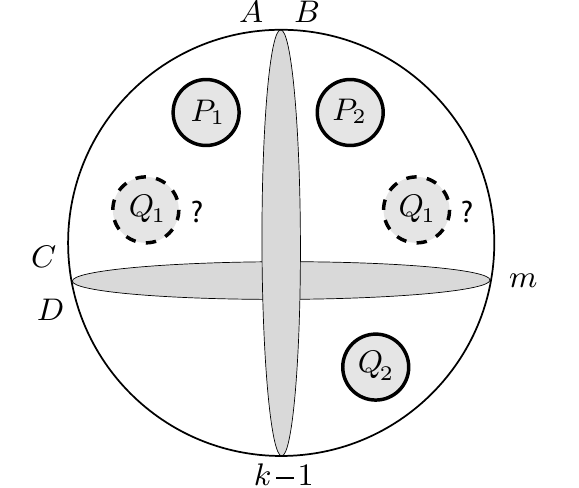}
   	  \caption{The known positions of $P_1, P_2, Q_1$ and~$Q_2$}
   \label{OptNestFig}
   \end{figure}

Now if $\XY := (A \cup C, B \cap D)$ has order~$<m$, then $\XY \in Q_2$ by~(P), and $\YX \in Q_1$ by consistency. Hence \XY\ distinguishes $Q_1$ from~$Q_2$ and has smaller order than \CD, contradicting the fact that \CD\ distinguishes $Q_1$ and~$Q_2$ efficiently. Thus \XY\ has order at least~$m$.

Hence by \eqref{ABorderCD} and Lemma~\ref{counting}, the order of $\EF := (B\cup D, A\cap C)$ is at most $k-1$. Then $\EF \in P_1$ by~(P), and $\FE \in P_2$ by consistency. Thus \EF\ distinguishes $P_1$ from~$P_2$. By \cite[Lemma~2.2]{confing},%
   \footnote{Swap the names of \CD\ and \EF\ in the statement of the lemma in~\cite{confing}.}
  \EF\ is nested with every separation in~\cN\ that \AB\ is nested with,%
   \COMMENT{}
   and in addition \EF\ is also nested with~\CD. Hence, \EF\ is nested with more separations in~\cN\ than \AB~is, contradicting the choice of~\AB.
\end{proof}

When we apply Lemma~\ref{lem_sopo_feasible} inductively, we have to make sure that every $\cN_\ell$ we construct consists only of $\cP_\ell$-essential separations. To ensure this, we have to reduce any task we tackle in the process of constructing~$\cN_\ell$. Given $k\ge 1$, a {\em $k$-strategy\/} is a $k$-tuple $(\sigma_1,\dots,\sigma_k)$ of strategies $\sigma_i$ each with range $\{\rops\}$. The restriction in the range of $k$-strategies will reduce our freedom in shaping the decompositions, but Example~\ref{the_Example} shows that considerable diversity remains.

Given $G$ and~$\cP$, a $k$-strategy $\Sigma=(\sigma_1,\dots,\sigma_k)$ \emph{determines\/} the set $\cN_\Sigma = \cN_\Sigma(G,\cP)$ defined recursively as follows. For $k=1$, let $\cN_\Sigma := \cN_{\sigma_1}\SP$, where $\cS$ is the set of proper $(<1)$-separations of~$G$. Then for $k\ge 2$ let 
\begin{equation}\label{eq_NSigma}
\cN_\Sigma := \cN\> \cup\!\! \bigcup_{O \in \cO_\cN}\! \cN_{\sigma_k}\SPi O\,,
\end{equation}
where $\cN = \cN_{\Sigma'}(G,\cP_{k-1})$ for $\Sigma' = (\sigma_1,\dots,\sigma_{k-1})$, and \cS\ is the set of proper $(k-1)$-separations of~$G$ that are nested with~\cN. The pairs $\SPi O$ are defined with reference to \cN, \cS\ and~$\cP = \cP_k$ as before Lemma~\ref{focus}.

As before, the sets $\cN_\Sigma$ will be {\em canonical\/} in that, for each~$\Sigma$, the map $(G,\cP) \mapsto \cN_\Sigma$ commutes with all isomorphisms $G\mapsto G'$.%
   \COMMENT{}
   In particular, if $\cP$ is invariant under the automorphisms of~$G$, then so is~$\cN_\Sigma$.

\begin{thm}\label{thm_can_sol_kP}
Every $k$-strategy $\Sigma$ determines for every set~\cP\ of $k$-profiles in a graph~$G$ a canonical nested system~$\cN_\Sigma(G,\cP)$ of $\cP$-essential separations of order~$<k$ that distinguishes all the profiles in~\cP\ efficiently.
\end{thm}
\begin{proof}
We show by induction on~$k$ that the recursive definition of $\cN_\Sigma(G,\cP)$ succeeds and that $\cN_\Sigma = \cN_\Sigma(G,\cP)$ has the desired properties. For $k=1$ this follows from Corollary~\ref{cor_kcon}.%
   \COMMENT{}

For $k \ge 2$ let \cN\ and \cS\ be defined as before the theorem. By the induction hypothesis, \cN\ is a nested system of $\cP_{k-1}$-essential separations of~$G$ that distinguishes the profiles in~$\cP_{k-1}$ efficiently. For every consistent orientation~$O$ of~\cN\ the pair \SPi O\ is a feasible task, by Lemma~\ref{lem_sopo_feasible}. By Theorem~\ref{thm_can_sol}, then, $\sigma_k$ determines a nested \sys\ $\cN_{\sigma_k}\SPi O \sub \cS_O \sub \cS$ that distinguishes all the profiles in~$\cP_O$. By definition of~\cS, all these $\cN_{\sigma_k}\SPi O$ are nested with~$\cN$, and they are nested with each other by Lemma~\ref{splitlemma}. Hence $\cN_\Sigma$ is well defined by~\eqref{eq_NSigma} and forms a nested \sys.

Let us show that $\cN_\Sigma$ has the desired properties. To show that $\cN_\Sigma$ distinguishes the profiles in~\cP\ efficiently, consider distinct $P,Q \in \cP$. If $\kappa(P,Q) < k-1$, then $P_{k-1} \neq Q_{k-1}$ are distinct profiles in $\cP_{k-1}$. So by the induction hypothesis there is a separation in $\cN \sub \cN_\Sigma$ that distinguishes $P$ from~$Q$ efficiently. If $\kappa(P,Q) = k-1$, we have $P_{k-1} = Q_{k-1}$. Then $P$ and~$Q$ have the same \cN-profile~$O$, and $P,Q \in \cP_O$. Hence there is a separation in $\cN_{\sigma_k}\SPi O \sub \cN_\Sigma$ that distinguishes $P$ from~$Q$; as it has order $k-1$, it does so efficiently.

It remains to show that every separation $\AB \in \cN_\Sigma$ is \cP-essential. If $\AB\in\cN$, this holds by the induction hypothesis and the definition of~$\cP_{k-1}$.%
   \COMMENT{}
   So assume that $\AB \in \cN_\Sigma \sm \cN$. Then there is a consistent orientation~$O$ of~\cN\ such that $\AB \in \cN_{\sigma_k}\SPi O$. Since $\sigma_k(i) \in \{\rops\}$ for all $i \in \N$, we know that \AB\ distinguishes some $P,Q\in\cP_O$. Then $\kappa(P,Q) = k-1 = |A\cap B|$, as otherwise \cN\ would distinguish $P$ from~$Q$ by the induction hypothesis. Hence, \AB\ distinguishes $P$ from $Q$ efficiently, as desired.
\end{proof}

It remains to translate our results from separation systems to \td s.
Recall from Theorem~\ref{treedec} that every nested proper \sys~\cN\ of~$G$ is induced by some \td~\TV: the separations of $G$ that correspond to edges of~\cT\ are precisely those in~\cN. In~\cite{confing} we showed that \TV\ is uniquely determined by~\cN.%
   \footnote{We assume here that parts corresponding to different nodes of~\cT\ are distinct. It is always possible to artificially enlarge the tree without changing the set of separations by duplicating a part.}
   Hence if~\cN\ is determined by some $k$-strategy, as in Theorem~\ref{thm_can_sol_kP}, we may say that this $k$-strategy {\em defines\/}~\TV\ on~$G$.%
   \COMMENT{}

If \cN\ comes from an application of Theorem~\ref{thm_can_sol_kP}, it will be proper (since it consists of essential separations~\eqref{distproper}) and canonical. In particular, if the set \cP\ of profiles considered is invariant under the automorphisms of~$G$, then so is~\cN, and hence so is~\cT: the automorphisms of $G$ will act on $V(\cT)$ as automorphisms of~\cT. And many natural choices of~\cP\ are invariant under the automorphisms of~$G$: the set of all $k$-profiles for given~$k$, the set of all $k$-block profiles, or the set of all tangles of order~$k$ to name some examples. All these can thus be distinguished in a single \td:

\begin{thm}\label{final}
Given $k\in\N$ and a graph~$G$, every $k$-strategy defines a canonical \td\ of adhesion $<k$ of~$G$ that distinguishes all its $k$-blocks and tangles of order~$k$. In particular, such a decomposition exists.\qed
\end{thm}

Theorem~\ref{final} is not the end of this story, but rather a beginning. One can now build on the fact that these \td s are given constructively and study their details. For example, we may wish to find out more about the structure or size of their parts, or obtain bounds on the number of parts containing a $k$-block or accommodating a tangle of order~$k$, compared with the total number of parts. The answers to such questions will depend on the $k$-strategy chosen. We shall pursue such questions in~\cite{CDHH13CanonicalParts}.

\bibliographystyle{plain}
\bibliography{collective}

\end{document}